\documentclass{article}
\usepackage{longtable}
\newcommand\tstrut{\rule{0pt}{4.4ex}}
\newcommand\bstrut{\rule[-3.0ex]{0pt}{0pt}}
\usepackage{caption}
\usepackage{latexsym,amsfonts,amssymb}
\usepackage[latin1]{inputenc}

\usepackage{tikz}
\usepackage{enumitem}
\usepackage{amsmath}
\usepackage{colortbl}
\usepackage{multirow}
\usepackage{lscape}
\usepackage{placeins}
\title{\sc Landau's Theorem on conjugacy classes for normal subgroups}

\author{Antonio Beltr\'an\\
\scriptsize
Departamento de Matem\'aticas,\\
\scriptsize Universidad Jaume I, \scriptsize
12071 Castell\'on, Spain\\
\scriptsize
e-mail: abeltran@mat.uji.es\\
\\
Mar\'{\i}a Jos\'e Felipe\\
\scriptsize
Instituto Universitario de Matem\'atica Pura y Aplicada,\\
\scriptsize Universidad Polit\'ecnica de Valencia, \scriptsize
46022 Valencia, Spain\\
\scriptsize e-mail: mfelipe@mat.upv.es \\
\\ Carmen Melchor\\
\scriptsize
Departamento de Matem\'aticas,\\
\scriptsize Universidad Jaume I, \scriptsize
12071 Castell\'on, Spain\\
\scriptsize
e-mail: cmelchor@uji.es      }

\date{}
\begin{document} \maketitle

\begin{abstract} Landau's theorem on conjugacy classes asserts that there are only finitely many finite groups, up to isomorphism, with exactly $k$ conjugacy classes for any positive integer $k$. We show that, for any positive integers $n$ and $s$, there exists only a finite number of finite groups $G$, up to isomorphism, having a normal subgroup $N$ of index $n$  which contains exactly $s$ non-central $G$-conjugacy classes. We provide upper bounds for the orders of $G$ and $N$, which are used by using GAP to classify all finite groups with normal subgroups having a small index and few $G$-classes. We also study the corresponding problems when we only take into account the set of $G$-classes of prime-power order elements contained in a normal subgroup.



\end{abstract}

\bigskip

\section{Introduction}
The renown Landau's theorem (\cite{Landau}) claims that there are only finitely many finite groups, up to isomorphism, with exactly $k$ conjugacy classes, for any positive integer $k$. However, no upper bound in terms of $k$ for the order of such groups was provided in the original proof  nor until 1963 by Brauer. A few years later, M. Newman (\cite{Newman}) demonstrated that $|G|\leq k^{2^{k-1}}$, or equivalently, that $k$ rises to infinity with the order of $G$ according to the inequality

$$k \geq \frac{{\rm log (} {\rm log } |G|)}{{\rm log } 4}.$$

Since then many results have appeared in the literature regarding logarithmic bounds for some special classes of groups. For instance, nilpotent or ``almost nilpotent" (\cite{CartwrightII}), or even solvable  groups (\cite{Cartwright}), as well as certain extensions of Landau's theorem, such as \cite{Kulshammer}, where only conjugacy classes of prime-power-order elements are taken into account. On the other hand,  a generalization of Landau's theorem for irreducible characters is proved in \cite{Lewis} by using the above result on classes of elements of prime-power order. \\

Let $G$ be a finite group and let $N$ be a normal subgroup of $G$. For each element $x \in N$, the $G$-conjugacy class is the set $x^{G}=\lbrace x^{g}\, |\, g \in G \rbrace$.  Bearing in mind that $N$ is union of $G$-classes, the aim of this paper is to extend Laudau's theorem for conjugacy classes of $N$. The natural  question is whether there exist finitely many groups having a normal subgroup which is union of a fixed number of $G$-classes or not. The answer is negative. In fact, if we take $N$ a $p$-elementary abelian group  of order $p^s$ and $G$ the holomorph group of $N$, then as a result of the fact that Aut$(N)$ acts transitively on $N\setminus\{1\}$, it follows that $N$ consists only of two $G$-classes, $\{1\}$ and $N\setminus\{1\}$. Nevertheless,  Aut$(N)\cong$ GL$(s, p)$ and so $|G :N|=|$GL$(s,p)|$ may increase as much as we want. In Section 2, we prove the following result which establishes that if we fix the index $|G:N|$, then the answer to our question is affirmative. We actually provide upper bounds for $|G|$ and $|N|$ depending on the number of non-central $G$-classes lying in $N$ (instead of all $G$-classes).\\

{\bf Theorem A}.
{\it Let $s, n\in \mathbb{N}$ such that $s,n \geq 1$. Then there exists at most a finite number of isomorphism classes of finite groups $G$ which contains a normal subgroup $N$ such that $|G:N|=n$ and $N$ has exactly $s$ non-central $G$-classes. Moreover,

 $$|G| < n^{2^s+1}(s+1)\prod_{i=0}^{s-1}(s+1-i)^{2^{s-1-i}}.$$
 and
  $$|N| < n^{2^s}(s+1)\prod_{i=0}^{s-1}(s+1-i)^{2^{s-1-i}}.$$
}\\
In the above theorem we take $s\geq 1$, since considering $s=0$ is equivalent to the fact that $N$ is central in $G$ and thus, a classification of such normal subgroups has no sense because these are infinitely many. When $n=1$, we obtain an improvement of Newman's bound in terms of the number of non-central classes in a group. \\

In Section 3, we present two applications of the above bounds which consist in classifying groups with normal subgroups having few non-central $G$-conjugacy classes. Particularly, we do this for normal subgroups having either one or two non-central $G$-classes of coprime size and for certain small indices by using GAP \cite{gap}. We recall that the classification of finite groups attending to their number of conjugacy classes has been a widely studied problem (see for instance \cite{Vera1} and \cite{Vera2}).\\

In Section 4, we study the problem of extending Landau's theorem for prime-power order elements lying in a normal subgroup $N$ of a group $G$. We prove that there also exist finitely many  finite groups, up to isomorphism, having a normal subgroup with fixed index and with  a fixed number of $G$-classes of prime-power order elements.   Really, we will see with examples that it is necessary to consider all $G$-classes instead of only non-central $G$-classes. Furthermore, we obtain a bound for the orders of $G$ and $N$ when $G$ is solvable.\\

{\bf Theorem B}. {\it For any two positive integers $k$ and $n$, there are only finitely many finite groups $G$, up to isomorphism, with a normal subgroup $N$ of index $n$ and containing $k$ $G$-classes of prime-power order elements. Moreover, when $G$ is solvable we have  $$|G|\leqslant n\prod_{i=0}^{nk-1}(nk-i)^{2^{i}} \quad {\it and} \quad  |N|\leqslant \prod_{i=0}^{nk-1}(nk-i)^{2^{i}}.$$}

\section{Proof of Theorem A}

We need the following elementary result, which is a slight improvement of a lemma in \cite{Newman}.\\

{\bf Lemma 2.1}. {\it Let $s, n\in \mathbb{N}$ with $s,n \geq 1$ and let $n_{1}, \cdots, n_{s+1}$ be integers such that $n_{1}\leq n_{2}\leq \cdots \leq n_{s+1}$ and $1/n=1/n_{1}+\cdots + 1/n_{s+1}$}. Then

$$n_{1}\leq n(s+1)$$ and
\begin{equation}\label{nk}
n_{k}\leq n^{2^{k-1}}(s+2-k)\prod_{i=0}^{k-2}(s+1-i)^{2^{k-2-i}}
\end{equation}
 for all $2\leq k\leq s+1$.\\

{\it Proof}. We obtain the first inequality trivially. We argue by induction on $k$ for $2\leq k\leq s+1 $. First, since $$1/n-1/n_{1}=1/n_{2}+\cdots+1/n_{s+1} \leq s/n_{2}$$ we deduce that $n_{2}\leq snn_{1}/(n_1-n)$ and we know that this denominator is bigger than 0 because $n_{1}>n$. By using the hypothesis we obtain $$n_{2}\leq snn_{1}$$ and since $n_{1}\leq n(s+1)$ we have  $$n_{2}\leq n^{2}(s+1)s.$$ Suppose that Eq. (\ref{nk}) holds for $k$ with $2\leq k \leq s$ and let us prove it for $k+1$. We know that $$0< 1/n-(1/n_{1}+\cdots + 1/n_{k})=1/n_{k+1}+\cdots +1/n_{s+1}\leq (s+1-k)/n_{k+1},$$ so we conclude that $$n_{k+1}\leq (s+1-k)nn_{1}\cdots n_{k}/(n_{1}\cdots n_{k}- (nn_2\cdots n_k+ \cdots +nn_{1}\cdots n_{k-1}))$$ and we know that this denominator is bigger than 0 because $1/n-(1/n_{1}+\cdots +1/n_{k})>0$. Thus, $$n_{k+1}\leq (s+1-k)nn_{1}\cdots n_{k}\leq$$ $$(s+1-k)nn(s+1)n^{2}(s+1)s\cdots n^{2^{k-1}}(s+2-k)\prod_{i=0}^{k-2}(s+1-i)^{2^{k-2-i}}.$$ Rearranging the terms on the right of the inequality it is easily computed that $$n_{k+1}\leq n^{2^{k}}(s+2-(k+1))\prod_{i=0}^{(k+1)-2}(s+1-i)^{2^{((k+1)-2)-i}}$$ and Eq.(1) is proved. $\Box$\\

{\it Proof of Theorem A.} Let $C_{i}=x_{i}^{G}$ for $i=1,\cdots, s$ the non-central $G$-conjugacy classes contained in $N$. Then $|C_{i}|=|G:{\rm \textbf{C}}_{G}(x_{i})|$. By the class equation we write $$|N|=|{\rm \textbf{Z}}(G)\cap N|+|C_{1}|+\cdots+|C_{s}|$$ and by dividing by $|G|$, we have $$1/n=1/|G:{\rm \textbf{Z}}(G)\cap N|+ 1/|{\rm \textbf{C}}_{G}(x_{1})|+\cdots 1/|{\rm \textbf{C}}_{G}(x_{s})|.$$ We rewrite the previous expression in the following way $$1/n=1/n_{1}+\cdots+1/n_{s+1}$$ where each $n_{i}$ is either $|{\rm\textbf{C}}_{G}(x_{j})|$ or $|G:{\rm \textbf{Z}}(G)\cap N|$. Without loss of generality, we can suppose that $n_{1}\leq \cdots \leq n_{s+1}$, so $$1/n=1/n_{1}+\cdots+1/n_{s+1} \leq (s+1)/n_{1}.$$ Thus,
\begin{equation}\label{n1}
n_{1}\leq n(s+1).
\end{equation}

We know that there exists some $i$ with $1\leq i\leq s+1$ which satisfies that $n_{i}=$ $|G:{\rm \textbf{Z}}(G)\cap N|$. Furthermore, for every $j\neq i$ we know that $n_{j}>|{\rm\textbf{Z}}(G)\cap N|$. Then $n_{i}> |G|/n_{j}$. Hence, $|G|< n_{i}n_{j}$ for every $j\neq i$, with $1\leq j\leq s+1$. In particular, if $i=1$, then $|G|< n_{1}n_{s+1}$ and if $i>1$ we have $|G|< n_{1}n_{i}\leq n_{1}n_{s+1}$. In all cases, by using Eq. (\ref{n1}) and by replacing $s+1$ in  Eq. (\ref{nk}) of Lemma 1, we get $$|G|< n_{1}n_{s+1}\leq n(s+1)n^{2^s}\prod_{i=0}^{s-1}(s+1-i)^{2^{s-1-i}}= n^{2^s+1}(s+1)\prod_{i=0}^{s-1}(s+1-i)^{2^{s-1-i}}. $$ The bound for $|N|$ trivially follows just by dividing by $n$.  Once the bound for $|G|$ is obtained, the first assertion in the Theorem is clear. $\Box$

\section{Applications of Landau's extension}
As an application of Theorem A we classify groups and normal subgroups with certain number of non-central classes of $G$ by using GAP. In order to expedite this classification, we use properties of the graph $\Gamma_{G}(N)$ associated to a normal subgroup $N$ of $G$ defined in \cite{Nuestro1}, as well as some results of \cite{Nuestro2} concerning to this graph when it has few vertices. These properties  give us relevant information about the structure of the normal subgroup. We recall that the graph $\Gamma_{G}(N)$ is defined in the following way: the  vertices are the non-central $G$-classes contained in $N$ and two vertices $x^{G}$ and $y^{G}$ are joined by an edge if and only if $|x^{G}|$ and $|y^{G}|$ have a common prime divisor. \\

We will show lists with the classification of those groups having a normal subgroup $N$ with one or two non-central $G$-conjugacy classes of coprime sizes for some concrete indices. The reason why we do not deal with the case of two non-central $G$-classes with non-coprime size is because we cannot assure that ${\rm \textbf{Z}}(G)\cap N=1$, whereas when the sizes are coprime we will see that this equality always holds, allowing us to improve the algorithm efficiency. \\

 The particular indices that we choose give rise to groups whose orders are included in the library SmallGroups of GAP. We provide information about the classified groups through tables in which appear: the index of the normal subgroup; the number of the corresponding SmallGroup and its group structure; the expression of the normal subgroup by using generators and its description; and finally, the sizes of the corresponding $G$-conjugacy classes.

\subsection{Normal subgroups with only one non-central $G$-class}
First, we show groups having a normal subgroup with only one non-central $G$-conjugacy class, what means that $\Gamma_{G}(N)$ has exactly one vertex. It is trivial that this case cannot happen for ordinary classes of $G$. \\

In pursuance of improving the algorithm efficiency, we particularize the bound of Theorem A and we use Theorem 3.1.2, which appeared in \cite{Nuestro2}. \\

{\bf Theorem 3.1.1.} {\it Let $N$ be a normal subgroup of a group $G$ with $|G:N|=n$. Suppose that $G$ has exactly only one non-central $G$-conjugacy class. Then $|G|<n(n+1)^2$.}\\

{\it Proof.} Let $x^G$ the only non-central $G$ conjugacy class of $N$. By the class equation we have

\begin{equation}\label{1class}
|N|=|N\cap \textbf{Z}(G)|+|x^{G}|.
\end{equation}

By dividing Eq. (\ref{1class}) by $|G|$ we have

 \begin{equation}\label{class1}
 \frac{1}{n}=\frac{1}{n_{1}}+\frac{1}{n_{2}}\leq \frac{2}{n_{1}}
 \end{equation}
 where $n_{i}$ is either $|G:N\cap \textbf{Z}(G)|$ or $|{\rm \textbf{C}}_{G}(x)|$, and we suppose without loss of generality that $n_{1}\leq n_{2}$. It follows that $n< n_{1}\leq 2n$. We observe that $n_{2}=nn_{1}/(n_{1}-n)$ is univocally determined by $n_{1}$. By arguing similarly as in the proof of Theorem A we get $|G|<n_{1}n_{2}=nn_{1}^{2}/(n_{1}-n)$, which is a continuous decreasing function for $n_{1}\in [n+1, 2n]$. We conclude that $|G|< n(n+1)^2$. $\Box$\\

In order to obtain our classification, we are able to discard those normal subgroups that do not satisfy the following conditions.\\

{\bf Theorem 3.1.2}.
{\it If $G$ is a group and $N$ is a normal subgroup of $G$ such that $\Gamma_{G}(N)$ has only one vertex, then $N$ is a $p$-group for some prime $p$ and $N/(N\cap{\bf Z}(G))$ is an elementary abelian $p$-group.}\\

{\it Proof.} See Theorem 3.1 of \cite{Nuestro2}. $\Box$ \\

In Table 1, we indicate for each index the bound for $|G|$ obtained in Theorem A, the corresponding improved bound of Theorem 3.1.1 and the number of groups with a normal subgroup containing a single non-central conjugacy class of $G$. In Table 2, we show the complete classification for the indices that appear in Table 1.\\

\FloatBarrier
\begin{table}[h]
\caption{Number of groups with normal subgroups having one non-central $G$-class.}
\centering
\begin{tabular}{ c  c  c c }
\hline
Index of $N$ & \parbox{2cm}{$|G|\leq 4n^3$ (Theorem A)} & \parbox{3cm}{$|G|\leq n(n+1)^2$ (Theorem 3.1.1)} & Number of groups \\
\hline \hline
2 & 32 & 18 & 3\\
3 & 108 & 48 & 2\\
4 & 256 & 100 & 21\\
5 & 500 & 180 & 0\\
6 & 864 & 294 & 16 \\
7 & 1372 & 448 & 1 \\
\hline
\end{tabular}
\end{table}
\FloatBarrier

The idea of the algorithm implemented in GAP is to list only integers less than or equal to the bound obtained in Theorem 3.1.1, say $c$, which are divisible by $n$, $n_{1}$ and $n_{2}$, appearing in that theorem, with $c/n_{1}>1$ and $c/n_{2}>1$. For $n=4$, among the possible orders that satisfy the above conditions, the order $64$ arises. This is the most conflictive possibility because there exist 267 groups and working with so many groups and normal subgroups is computationally expensive. However, we can discard this case by proving the following property generalized for $p$-groups. For $n=8$, in some cases, specifically $2$-groups, although the order is included in the library, we have had to stop the classification because the amount of possibilities is too large.\\

{\bf Theorem 3.1.3}.
{\it Let $G$ be a group such that $|G|=p^k$ with $p$ prime, and let $N$ be a normal subgroup of $G$ of order $p^a$ such that $N$ has only one non-central $G$-conjugacy class. Then $p=2$, $N$ is abelian and $a\leq (k+1)/2$.}\\

{\it Proof.} Suppose that $G$ has a normal subgroup with one single non-central $G$-conjugacy class, say $|x^G|$. It follows that $|N|=|{\rm \textbf{Z}}(G)\cap N |+|x^G|$ what means that $p^a=p^c+p^d=p^c(1+p^{d-c})$ with $c\leq d$. This equation only holds if $p=2$. Accordingly, we obtain $2^{a}=2^c(1+2^{d-c})$ with $c\leq d$. This forces to $c=d=a-1$, so $|{\rm \textbf{Z}}(G)\cap N|=|x^G|=2^{a-1}$. Now, since ${\rm \textbf{Z}}(G)\cap N \subseteq {\rm \textbf{Z}}(N)$ and $|N:{\rm \textbf{Z}}(G)\cap N|=2$ we conclude that $N$ is abelian. Then, $N\subseteq {\rm \textbf{C}}_{G}(x)$ and, since  $2^{k-1}=|x^G|=|G:{\rm \textbf{C}}_G(x)|$ divides $|G:N|=2^{k-a}$, we obtain $a-1\leq k-a$, that is, $a\leq (k+1)/2$. $\Box$

\FloatBarrier
\begin{small}
\begin{center}
\begin{longtable}[t]{ l  l  l l}
\caption{Normal subgroups with one non-central $G$-class for the indexes of Table 1.}\\
\hline
$n$ & $G$ & $N$ & $|x^{G}|$  \\
\hline
\textbf{2} & \parbox{3cm}{$SmallGroup(6, 1)=S_{3}$} &$C_{3}\cong Group([ f2 ])$ & 2 \tstrut\bstrut \\
\cline{2-4}
& \multirow{3}{*}{\parbox{3cm}{$SmallGroup(8, 3)=D_{8}$}}  &$C_{4}\cong Group([ f1*f2, f3 ])$ & 2\tstrut\\
& &$C_{2}\times C_{2}\cong Group([ f1, f3 ])$ & 2   \\
& &$C_{2}\times C_{2}\cong Group([ f2, f3 ])$ & 2   \bstrut\\
\cline{2-4}
& \multirow{3}{*}{\parbox{3cm}{$SmallGroup(8, 4)=Q_{8}$}}  & $C_{4}\cong Group([ f1*f2, f3 ])$ &  2\tstrut\\
& &  $C_{4}\cong Group([ f1, f3 ])$ &  2\\
& &  $C_{4}\cong Group([ f2, f3 ])$ &  2\bstrut\\
\hline
\textbf{3} & \parbox{3cm}{$SmallGroup(12, 3)=A_{4}$} & $C_{2}\times C_{2}\cong Group([ f2, f3 ])$ & 3\tstrut\bstrut\\
\cline{2-4}
& \parbox{3cm}{$SmallGroup(24, 3)=SL(2,3)$} & $Q_{8}\cong Group([ f2, f3, f4 ])$ & 6 \tstrut\bstrut\\
\hline
\textbf{4} & \parbox{3cm}{$SmallGroup(12, 1)=C_3 : C_4$} & $C_3\cong Group([ f3 ])$ & 2 \tstrut\bstrut\\
\cline{2-4}
& \parbox{3cm}{$SmallGroup(12, 4)=D_{12}$} & $C_3\cong Group([ f3 ]), $ & 2  \tstrut\bstrut\\
\cline{2-4}
& \multirow{2}{*}{\parbox{3.5cm}{$SmallGroup(16, 3)=(C_4 \times C_2) : C_2$}} & $C_2 \times C_2\cong Group([ f2*f4, f3 ])$ & 2\tstrut\\
& & $C_2 \times C_2\cong Group([ f2, f3 ])$ & 2\bstrut\\
\cline{2-4}
& \multirow{2}{*}{\parbox{3cm}{$SmallGroup(16, 4)=C_4 : C_4$}} & $C_4\cong Group([ f2*f4, f3 ])$ & 2\tstrut\\
& & $C_4\cong Group([ f2, f3 ]), $ & 2 \bstrut\\
\cline{2-4}
& \multirow{2}{*}{\parbox{3cm}{$SmallGroup(16, 6)=C_8 : C_2$}} & $C_4\cong Group([ f2*f3, f4 ])$ & 2\tstrut\\
& & $C_2\times C_2\cong Group([ f2, f4 ])$ & 2 \bstrut\\
\cline{2-4}
& \parbox{3cm}{$SmallGroup(16, 7)=D_{16}$} & $C_4\cong Group([ f3, f4 ])$ & 2 \tstrut\bstrut\\
\cline{2-4}
& \parbox{3cm}{$SmallGroup(16, 8)=QD_{16}$} & $C_4\cong Group([ f3, f4 ])$ & 2 \tstrut\bstrut\\
\cline{2-4}
& \parbox{3cm}{$SmallGroup(16, 9)=Q_{16}$} & $C_4\cong Group([ f3, f4 ])$ & 2 \tstrut\bstrut\\
\cline{2-4}
& \multirow{6}{*}{\parbox{3cm}{$SmallGroup(16, 11)=C_{2}\times D_{8}$}} & $C_4\cong Group([ f1*f2*f3, f4 ])$ &2  \tstrut\\
& & $C_4\cong Group([ f1*f2, f4 ])$ & 2   \\
& & $C_2\times C_2\cong Group([ f1*f3, f4 ])$ & 2   \\
& & $C_2\times C_2\cong Group([ f1, f4 ])$ & 2   \\
& & $C_2\times C_2\cong Group([ f2*f3, f4 ])$ & 2   \\
& & $C_2\times C_2\cong Group([ f2, f4 ])$ & 2   \bstrut\\
\cline{2-4}
& \multirow{6}{*}{\parbox{3cm}{$SmallGroup(16, 12)=C_{2}\times Q_{8}$}} & $C_4\cong Group([ f1*f2*f3, f4 ])$ & 2   \tstrut\\
&  & $C_4\cong Group([ f1*f2, f4 ])$ & 2   \\
&  & $C_4\cong Group([ f1*f3, f4 ])$ & 2   \\
&  & $C_4\cong Group([ f1, f4 ])$ & 2   \\
&  & $C_4\cong Group([ f2*f3, f4 ])$ & 2   \\
&  & $C_4\cong Group([ f2, f4 ])$ & 2  \bstrut\\
\cline{2-4}
& \multirow{6}{*}{\parbox{3cm} {$SmallGroup(16, 13)=(C_{4}\times C_{2}):C_{2}$}} & $C_2\times C_2\cong Group( [f1*f2*f3, f4 ])$ & 2    \tstrut\\
& & $C_4\cong Group([ f1*f2, f4 ])$ & 2   \\
& & $C_4\cong Group([ f1*f3, f4 ])$ & 2   \\
& & $C_2\times C_2=Group([ f1, f4 ])$ & 2   \\
& & $C_4\cong Group([ f2*f3, f4 ])$ & 2   \\
& & $C_2\times C_2\cong Group([ f2, f4 ])$ & 2   \bstrut\\
\cline{2-4}
& \parbox{3cm}{$SmallGroup(20, 3)=C_{5}: C_{4}$} & $C_5\cong Group([ f3 ])$ & 4  \tstrut\bstrut\\
\cline{2-4}
& \multirow{4}{*}{\parbox{3.5cm} {$SmallGroup(32, 27)=(C_2 \times C_2 \times C_2 \times C_2) : C_2$}} & $C_4 \times C_2\cong Group([ f1*f2*f3, f4, f5 ])$ & 4   \tstrut\\
&  & $C_4 \times C_2\cong Group([ f1*f2, f4, f5 ])$ & 4   \\
&  & $C_4 \times C_2\cong Group([ f1*f3, f4, f5 ])$ & 4   \\
&  & $C_2 \times C_{2} \times C_2\cong Group([ f1, f4, f5 ])$ & 4    \bstrut\\
\cline{2-4}
& \multirow{5}{*}{\parbox{3cm} {$SmallGroup(32, 28)=(C_4 \times C_2 \times C_2) : C_2$}} & $C_4 \times C_{2}\cong Group([ f1*f2*f3, f4, f5 ])$ & 4   \tstrut\\
&  & $C_2 \times C_{2}\times C_{2}\cong Group([ f1*f2, f4, f5 ])$ & 4   \\
&  & $C_4 \times C_{2}\cong Group([ f1*f3, f4, f5 ])$ & 4   \\
&  & $C_2 \times C_{2}\times C_{2}\cong Group([ f1, f4, f5 ])$ & 4   \bstrut\\
\cline{2-4}
& \multirow{4}{*}{\parbox{3cm} {$SmallGroup(32, 29)= (C_2 \times Q_8) : C_2$}} & $C_4 \times C_{2}\cong Group([ f1*f2*f3, f4, f5 ])$ & 4   \tstrut\\
&  & $C_4 \times C_{2}\cong Group([ f1*f2, f4, f5 ])$ & 4   \\
&  & $C_4 \times C_{2}\cong Group([ f1*f3, f4, f5 ])$ & 4   \\
&  & $C_4 \times C_{2}\cong Group([ f1, f4, f5 ])$ & 4    \bstrut\\
\cline{2-4}
& \multirow{4}{*}{\parbox{3cm} {$SmallGroup(32, 30) = (C_4 \times C_2 \times C_2) : C_2$}} & $C_4 \times C_{2}\cong Group([ f1*f2*f3, f4, f5 ])$ & 4   \tstrut\\
& & $C_4 \times C_{2}\cong Group([ f1*f2, f4, f5 ])$ & 4   \\
& & $C_4 \times C_{2}\cong Group([ f1*f3, f4, f5 ])$ & 4   \\
&  & $C_2 \times C_{2}\times C_{2}\cong Group([ f1, f4, f5 ])$ & 4    \bstrut\\
\cline{2-4}
& \multirow{5}{*}{\parbox{3cm} {$SmallGroup(32, 31) =(C_4 \times C_4) : C_2$}} & $C_2 \times C_{2}\times C_{2}\cong Group([ f1*f2*f3, f4, f5 ])$ & 4  \tstrut\\
& & $C_4 \times C_{2}\cong Group([ f1*f2, f4, f5 ])$ & 4   \\
& & $C_4 \times C_{2}\cong Group([ f1*f3, f4, f5 ])$ & 4   \\
& & $C_2 \times C_{2}\times C_{2}\cong Group([ f1, f4, f5 ])$ & 4    \bstrut\\
\cline{2-4}
& \multirow{4}{*}{\parbox{3.3cm} {$SmallGroup(32, 32) =(C_2 \times C_2) . (C_2 \times C_2 \times C_2)$}} & $C_4 \times C_{2}\cong Group([ f1*f2*f3, f4, f5 ])$ & 4   \tstrut\\
& & $C_4 \times C_{2}\cong Group([ f1*f2, f4, f5 ])$ & 4   \\
& & $C_4 \times C_{2}\cong Group([ f1*f3, f4, f5 ])$ & 4   \\
& & $C_4 \times C_{2}\cong Group([ f1, f4, f5 ])$ & 4   \bstrut\\
 \cline{2-4}
& \multirow{4}{*}{\parbox{3cm} {$SmallGroup(32, 33) =(C_4 \times C_4) : C_2$}} & $C_4 \times C_{2}\cong([ f1*f2*f3, f4, f5 ])$ & 4  \tstrut\\
& & $C_4 \times C_{2}\cong([ f1*f2, f4, f5 ])$ & 4   \\
& & $C_4 \times C_{2}\cong([ f1*f3, f4, f5 ])$ & 4   \\
& & $C_2 \times C_{2}\times C_{2}\cong([ f1, f4, f5 ])$ & 4    \bstrut\\
 \cline{2-4}
& \multirow{4}{*}{\parbox{3cm} {$SmallGroup(32, 34) =(C_4 \times C_4) : C_2$}} & $C_2 \times C_{2}\times C_{2}\cong([ f1*f2*f3, f4, f5 ])$ & 4  \tstrut\\
& & $C_2 \times C_{2}\times C_{2}\cong([ f1*f2, f4, f5 ])$ & 4   \\
& & $C_2 \times C_{2}\times C_{2}\cong([ f1*f3, f4, f5 ])$ & 4   \\
& & $C_2 \times C_{2}\times C_{2}\cong([ f1, f4, f5 ])$ & 4    \bstrut\\
 \cline{2-4}
& \multirow{4}{*}{\parbox{3cm} {$SmallGroup(32, 35) =(C_4 : Q_8) : C_2$}} & $C_4 \times C_{2}\cong([ f1*f2*f3, f4, f5 ])$ & 4  \tstrut\\
& & $C_4 \times C_{2}\cong([ f1*f2, f4, f5 ])$ & 4   \\
& & $C_4 \times C_{2}\cong([ f1*f3, f4, f5 ])$ & 4   \\
& & $C_4 \times C_{2}\cong([ f1,f4,f5 ])$ & 4    \bstrut\\
\hline
\textbf{6} & \parbox{3cm}{$SmallGroup(18, 1)=D_{18}$} &$C_{3}\cong Group([ f3 ])$ & 2 \tstrut\bstrut \\
\cline{2-4}
& \parbox{3cm}{$SmallGroup(18, 3)=C_{3}\times S_{3}$} &$C_{3}\cong Group([ f3 ])$ & 2 \tstrut\bstrut \\
\cline{2-4}
& \multirow{4}{*}{\parbox{3cm}{$SmallGroup(18, 4)=(C_{3} \times C_{3}) : C2$}}  &$C_{3}\cong Group([ f2*f3^2 ])$ & 2\tstrut\\
& &$C_{3}\cong Group([ f2*f3 ])$ & 2   \\
& &$C_{3}\cong Group([ f2 ])$ & 2   \\
& &$C_{3}\cong Group([ f3 ])$ & 2   \bstrut\\
\cline{2-4}
& \parbox{3cm}{$SmallGroup(24, 12)=S_{4}$} &$C_{2}\times C_{2}\cong Group([ f3, f4 ])$ & 3 \tstrut\bstrut \\
\cline{2-4}
& \parbox{3cm}{$SmallGroup(24, 13)=C_{2}\times A_{4}$} &$C_{2}\times C_{2}\cong Group([ f3, f4 ])$ & 3 \tstrut\bstrut \\
\cline{2-4}
& \parbox{3cm}{$SmallGroup(24, 4)=C_{3}: Q_{8}$} &$C_{4}\cong Group([ f2, f3 ])$ & 2 \tstrut\bstrut \\
\cline{2-4}
& \parbox{3cm}{$SmallGroup(24, 6)=D_{24}$} &$C_{4}\cong Group([ f2, f3 ])$ & 2 \tstrut\bstrut \\
\cline{2-4}
& \parbox{3cm}{$SmallGroup(24, 8)=(C_{6} \times C_{2}) : C_{2}$} &$C_{2}\times C_{2}\cong Group([ f2, f3 ])$ & 2 \tstrut\bstrut \\
\cline{2-4}
& \multirow{3}{*}{\parbox{3cm}{$SmallGroup(24, 10)=C_{3}\times D_{8}$}}  & $C_{4}\cong Group([ f1*f2, f4 ])$ &  2\tstrut\\
& &  $C_{2}\times C_{2}\cong Group([ f1, f4 ])$ &  2\\
& &  $C_{2}\times C_{2}\cong Group([ f2, f4 ])$ &  2\bstrut\\
\cline{2-4}
& \multirow{3}{*}{\parbox{3cm}{$SmallGroup(24, 11)=C_{3}\times Q_{8}$}}  & $C_{4}\cong Group([ f1*f2, f4 ])$ &  2\tstrut\\
& &  $C_{4}\cong Group([ f1, f4 ])$ &  2\\
& &  $C_{4}\cong Group([ f2, f4 ])$ &  2\bstrut\\
\cline{2-4}
& \parbox{3cm}{$SmallGroup(42, 1)=(C_{7} : C_{3}) : C_{2}$} &$C_{7}\cong Group([ f3 ])$ & 6 \tstrut\bstrut \\
\cline{2-4}
& \parbox{3cm}{$SmallGroup(48, 28)=SL(2,3) . C_{2}$} &$Q_{8}\cong Group([ f3, f4, f5 ])$ & 6 \tstrut\bstrut \\
\cline{2-4}
& \parbox{3cm}{$SmallGroup(48, 29)=GL(2,3)$} &$Q_{8}\cong Group([ f3, f4, f5 ])$ & 6 \tstrut\bstrut \\
\cline{2-4}
& \parbox{3cm}{$SmallGroup(48, 32)=C_{2} \times SL(2,3)$} &$Q_{8}\cong  Group([ f3, f4, f5 ])$ & 6 \tstrut\bstrut \\
\cline{2-4}
& \parbox{3cm}{$SmallGroup(48, 33)=SL(2,3) : C_{2}$} &$Q_{8}\cong Group([ f3, f4, f5 ])$ & 6 \tstrut\bstrut \\
\cline{2-4}
& \multirow{4}{*}{\parbox{3cm}{$SmallGroup(54, 8)=((C_{3} \times C_{3}) : C_{3}) : C_{2}$}}  & $C_{3}\times C_{3}\cong Group([ f2*f3^2, f4 ])$ &  6\tstrut\\
& &  $C_{3}\times C_{3}\cong  Group([ f2*f3, f4 ])$ &  6\\
& &  $C_{3}\times C_{3}\cong Group([ f2, f4 ])$ &  6\\
& &  $C_{3}\times C_{3}\cong Group([ f3, f4 ])$ &  6\bstrut\\
\cline{2-4}
\textbf{7} & \parbox{3cm}{$SmallGroup(56, 11)=(C_{2} \times C_{2} \times C_{2}) : C_{7}$} &$C_{2} \times C_{2} \times C_{2}\cong Group([ f2, f3, f4 ])$ & 7 \tstrut\bstrut \\
\hline
\end{longtable}
\end{center}
\end{small}
\FloatBarrier

{\bf Remark 1.} The only case in which $N$ is not abelian is when $N\cong Q_{8}$. This happens exactly for $n=3$ and $G=SL(2,3)$, and for $n=6$ and $4$ groups of order $48$. \\

\subsection{Normal subgroups with two non-central $G$-classes of coprime sizes}
In this section, we study normal subgroups with two non-central $G$-classes of coprime sizes. It was proved in \cite{BerHerMann} that this case for ordinary classes only happens when $G\cong S_{3}$. However, for $G$-classes lying in a normal subgroup this is not the unique possibility (see \cite{Nuestro2}). We use the following in order to improve the efficiency of the algorithm in this case.\\

{\bf Lemma 3.2.1.} {\it Let $N$ be a normal subgroup of a group $G$. If $\Gamma_{G}(N)$ has two vertices but no edges, then {\rm \textbf{Z}}$(G)\cap N=1$.}\\

{\it Proof}. See Lemma 4.1.1 of \cite{Nuestro2}. $\Box$\\

The above let us improve the bounds of Theorem A.\\

{\bf Theorem 3.2.2.} {\it Let $N$ be a normal subgroup of a group $G$ with $|G:N|=n$. Suppose that $G$ has exactly two non-central conjugacy classes $x_{1}^{G}$ and $x_{2}^{G}$ of coprime sizes lying in $N$ and let $n_{1}=|${\rm \textbf{C}}$_{G}(x_{1})|$ and $n_{2}=|${\rm \textbf{C}}$_{G}(x_{2})|$ such that $n_{1}<n_{2}$. Then,}

\begin{enumerate}
\item[(i)] {\it $n+1 \leq n_{1} \leq 3n -1$.}
\item[(ii)] {\it$E[\frac{nn_{1}}{n_{1}-n}] + 1\leq n_{2}\leq E[\frac{2n_{1}n}{n_{1}-n}]-1$, where $E(x)$ denotes the integer part of $x$.}
\item[(iii)] {\it$|G|\leq n(n+1)(n^2+n+1)$.}
\end{enumerate}

{\it Proof.} By the class equation and Lemma 3.2.1 we write
\begin{equation}\label{class}
|N|=1+|x_{1}^{G}|+|x_{2}^{G}|.
\end{equation}
By dividing Eq. (\ref{class}) by $|G|$,  we have

 \begin{equation}\label{class1}
 \frac{1}{n}=\frac{1}{|G|}+\frac{1}{n_{1}}+\frac{1}{n_{2}}<\frac{3}{n_{1}}.
 \end{equation}
Therefore, $n_{1}<3n$ and $n<n_{1}$, so (i) holds. Furthermore, by noticing that $n < n_{1} < n_{2} < |G|,$ we get

 \begin{equation*}
 \frac{1}{n}-\frac{1}{n_{1}}=\frac{1}{n_{2}}+\frac{1}{|G|}< \frac{2}{n_{2}}
 \end{equation*}
  and this means that
 \begin{equation*}
 n_{2}< E[\frac{2n_{1}n}{n_{1}-n}].
 \end{equation*}

On the other hand,

 \begin{equation*}
 \frac{1}{n}-\frac{1}{n_{1}}=\frac{1}{n_{2}}+\frac{1}{|G|} > \frac{1}{n_{2}},
 \end{equation*}
 so
 \begin{equation*}
 n_{2} > E[\frac{nn_{1}}{n_{1}-n}].
 \end{equation*} Observe that $\frac{nn_{1}}{n_{1}-n}=n+\frac{n^{2}}{n_{1}-n}>n>n_{1}$. Then (ii) is proved.\\

Also, Eq. (\ref{class1}) yields to
\begin{equation}\label{eq}
|G|=nn_{1}n_{2}/(n_{1}n_{2}-nn_{2}-nn_{1})
 \end{equation}
 Moreover, we know that $n_{1}n_{2}-nn_{2}-nn_{1}>0$ because $1/n-1/n_{1}-1/n_{2}>0$. Then Eq. (\ref{eq}) defines a continuous function when $n_{1}$ and $n_{2}$ belong to the intervals of (i) and (ii). It easy to check that this function reaches its maximum for the minimum values of $n_{1}$ and $n_{2}$. Thus, we obtain (iii) by taking the lower bounds for $n_{1}$ and $n_{2}$, that is, $n_{1}=n+1$ and $n_{2}=n(n+1)+1$. By replacing these in Eq. (\ref{eq}) we have $|G|\leq n(n+1)(n^2+n+1)$. $\Box$ \\

Moreover, by using the following result of \cite{Nuestro2} we can rule out, for each possible order, those normal subgroups that do not have the structure described in the following.\\

{\bf Theorem 3.2.3.} {\it Let $N$ be a normal subgroup of a group $G$ such that $\Gamma_G(N)$ has two vertices and no edge. Then $N$ is a $2$-group or a Frobenius group with  $p$-elementary abelian
kernel $K$, and complement $H$, which is cyclic of order $q$, for two different primes $p$ and $q$. In particular, $|N|=p^n q$ with $n\geq 1$.}\\

{\it Proof.} This is Theorem 4.1.2 of \cite{Nuestro2}. $\Box$\\

In Table 3 we indicate, for each index, the bound for $|G|$ obtained in Theorem A, an improved bound by using Theorem 3.2.2 and the number of groups with a normal subgroup containing two non-central conjugacy class of $G$ of coprime sizes. In Table 4, we show the complete classification for the indices that appear in Table 3. \\

 \FloatBarrier
\begin{table}[h]
\begin{small}
\caption{Groups with normal subgroups having two non-central $G$-classes with coprime sizes.}
\centering
\begin{tabular}{ c  c  c  c }
\hline
Index of $N$ &  \parbox{2cm}{$|G|\leq54n^5$ (Theorem A)} & \parbox{3.5cm}{$|G|\leq n(n+1)(n^2+n+1) $ (Theorem 3.2.2)} & Number of groups\\
\hline \hline2 & 1728 & 42 & 3\\
3 & 13122 & 156 & 2\\
4 & 55296 & 420 & 7\\
5 & 168750 & 930 & 2\\
6 & 419904 & 1806 & 8\\
7 & 907578 & 3192 & 1\\
8 & 1769472 & 5256 & 22\\
9 & 3188646 & 8190 & 5\\
10 & 5400000 & 12210 & 7\\
11 & 8696754 & 17556 & 2\\
13 & 20049822 & 33306 & 2\\
17 & 76672278 & 93942 & 1\\
\hline
\end{tabular}
\end{small}
\end{table}
\FloatBarrier

 The idea of the algorithm implemented in GAP is to list only those integers less than or equal to the bound obtained in Theorem 3.2.2, say $c$, satisfying that $n$,  $n_{1}$ and $n_{2}$ divide $c$, with $(c/n_{1}, c/n_{2})=1$, $c/n_{1}>1$ and $c/n_{2}>1$. For $n=12$,  we have not been able to give a complete classification because $2184$ is a possible value for $|G|$ that satisfies such conditions, and however, this order is not included in the SmallGroups library of GAP. The same problem happens for the indices $14$, $15$, $16$ and $18$.

\FloatBarrier
\begin{footnotesize}
\begin{center}
\begin{longtable}[h]{l l  l l l}
  \caption{Normal subgroups with two non-central $G$-classes with coprime sizes whose index is specified in Table 3.}\\
\hline
$n$ & $G$  & $N$ & $|x^{G}|$ & $|y^{G}|$  \\
\hline \hline
\textbf{2} & \parbox{3cm} {$SmallGroup(24, 12)=S_{4}$} & $A_{4} \cong Group([ f2, f3, f4 ])$ & 8 &3  \tstrut\bstrut \\
\cline{2-5}
& \parbox{3cm} {$SmallGroup(20, 3)=C_{5}:C_{4}$} & $D_{10}\cong Group([ f2, f3 ])$ & 5 & 4  \tstrut\bstrut \\
\cline{2-5}
& \multirow{2}{*}{\parbox{3cm} {$SmallGroup(12, 4)=D_{12}$}} & $S_{3}\cong Group([ f1*f2, f3 ])$ & 3 & 2 \tstrut\\
&  &  $S_{3}\cong Group([ f1, f3 ])$ & 3 & 2 \bstrut\\
\hline
\textbf{3} & \parbox{3cm} {$SmallGroup(42, 1)=(C_{7}:C_{3}):C_{2}$}  & $D_{14}\cong Group([ f1, f3 ])$ & 7 & 6\tstrut\bstrut \\
\cline{2-5}
& \parbox{3cm} {$SmallGroup(18, 3)=C_{3}\times S_{3}$} & $S_{3}\cong Group([ f1, f3 ])$ & 3 & 2  \tstrut\bstrut \\
\hline
\textbf{4} & \parbox{3cm} {$SmallGroup(48, 30)=A_{4}:C_{4}$} & $A_{4}\cong Group([ f3, f4, f5 ])$& 8 &3 \tstrut\bstrut \\
\cline{2-5}
& \parbox{3cm} {$SmallGroup(48, 48)=C_{2} \times S_{4}$} & $A_{4}\cong Group([ f3, f4, f5 ])$ & 8 &3  \tstrut\bstrut \\
\cline{2-5}
& \parbox{3cm} {$SmallGroup(72, 39)=(C_{3} \times C_{3}): C_{8}$} & $(C_{3} \times C_{3}): C_{2}\cong Group([ f3, f4, f5 ])$ & 9 &8\tstrut\bstrut \\
\cline{2-5}

& \parbox{3cm} {$SmallGroup(72, 41)=(C_{3} \times C_{3}): Q_{8}$} & $(C_{3} \times C_{3}): C_{2}\cong Group([ f3, f4, f5 ])$ & 9 &8 \tstrut\bstrut \\
\cline{2-5}

& \multirow{2}{*}{\parbox{3cm} {$SmallGroup(40, 12)=C_{2}\times (C_{5} : C_{4})$}} & $D_{10}\cong Group([ f2*f3, f4 ])$ & 5 &4 \tstrut \\
&  &  $D_{10}\cong Group([ f3, f4 ])$ & 3 & 2 \bstrut \\
\cline{2-5}

& \multirow{2}{*}{\parbox{3cm} {$SmallGroup(24, 5)=C_{4}\times S_{3}$}} & $S_{3}\cong Group([ f1*f3, f4 ])$ & 3 &2  \tstrut \\
&  & $S_{3}\cong Group([ f1,  f4 ])$ & 3 & 2 \bstrut \\
\cline{2-5}

& \multirow{4}{*}{\parbox{3cm} {$SmallGroup(24, 14)=C_{2}\times C_{2} \times S_{3}$}} & $S_{3}\cong Group([ f1*f2*f3, f4 ])$ & 3 &2  \tstrut \\
&  & $S_{3}\cong Group([ f1*f2, f4 ])$ & 3 & 2 \\
&  & $S_{3}\cong Group([f1*f3, f4 ])$ & 3 & 2  \\
&  & $S_{3}\cong Group([ f1, f4 ])$ & 3 & 2\bstrut \\
\hline
\textbf{5} & \parbox{3cm} {$SmallGroup(110, 1)=(C_{11}:C_{5}):C_{2}$} &  $D_{22}\cong Group([ f1, f3 ])$
 & 11 & 10 \tstrut\bstrut \\
\cline{2-5}
& \parbox{3cm} {$SmallGroup(30, 1)=C_{5}\times S_{3}$} & $S_{3}\cong Group([ f1, f3 ])$ & 3 & 2 \tstrut\bstrut \\
\hline
\textbf{6} & \parbox{3cm} {$SmallGroup(72, 42)=C_{3}\times S_{4}$} & $A_{4}\cong Group([ f3, f4, f5 ])$ & 8 & 3 \tstrut\bstrut \\
\cline{2-5}
& \multirow{3}{*}{\parbox{3cm} {$SmallGroup(72, 43)=(C_{3}\times A_{4}): C_{2}$}} & $A_{4}\cong Group([ f2*f3^2, f4, f5 ])$ & 8 & 3\tstrut\\
&  & $A_{4}\cong Group([ f2*f3, f4, f5 ])$ & 8 & 3  \\
&  & $A_{4}\cong Group([ f2, f4, f5 ])$ & 8 & 3\bstrut \\
\cline{2-5}
& \parbox{3cm} {$SmallGroup(156, 7)=(C_{13} : C_{4}) : C_{3}$} & $D_{26}\cong Group([ f3, f4 ])$ & 13 & 12\tstrut\bstrut \\
\cline{2-5}
& \multirow{2}{*}{\parbox{3cm} {$SmallGroup(84, 7)=C_{2}\times (C_{7}:C_{3}):C_{2}$}} & $D_{14}\cong Group([ f1*f2, f4 ])$ & 7 & 6\tstrut \\
&  & $D_{14}\cong Group([ f1, f4 ])$ & 7 & 6\bstrut \\
\cline{2-5}
& \parbox{3cm} {$SmallGroup(60, 6)=C_{3}\times (C_{5}: C_{4})$} & $D_{10}\cong Group([ f3, f4 ])$ & 5 & 4\tstrut\bstrut \\
\cline{2-5}
& \parbox{3cm} {$SmallGroup(60, 7)=C_{15} : C_{4}$} & $D_{10}\cong Group([ f2, f4 ])$ & 5 & 4\tstrut\bstrut \\
\cline{2-5}
& \multirow{2}{*}{\parbox{3cm} {$SmallGroup(36, 10)=S_{3} \times S_{3}$}} & $S_{3}\cong Group([ f2, f3 ])$ & 3 & 2\tstrut\\
&  & $S_{3}\cong Group([ f1, f4 ])$ & 3 & 2\bstrut \\
\cline{2-5}
& \multirow{2}{*}{\parbox{3cm} {$SmallGroup(36, 12)=C_{6} \times S_{3}$}} & $S_{3}\cong Group([ f1*f2, f4 ])$ & 3 & 2\tstrut\\
&  & $Group([ f1, f4 ]$ & 3 & 2\bstrut \\
\hline
\textbf{7} & \parbox{3cm} {$SmallGroup(42, 3)=C_{7}\times S_{3}$} & $S_{3}\cong Group([ f1, f3 ])$ & 3 & 2 \tstrut\bstrut \\
\hline
\textbf{8} & \parbox{3cm} {$SmallGroup(96, 65)=A_{4}: C_{8}$}& $A_{4}\cong Group([ f4, f5, f6 ])$ & 8 & 3 \tstrut\bstrut \\
\cline{2-5}

& \parbox{3cm} {$SmallGroup(96, 185)=A_{4}: Q_{8}$} & $A_{4}\cong Group([ f4, f5, f6 ])$ & 8 & 3\tstrut\bstrut \\
\cline{2-5}

& \parbox{3cm} {$SmallGroup(96, 186)=C_{4} \times S_{4}$} & $A_{4}\cong Group([ f4, f5, f6 ])$ & 8 & 3\tstrut\bstrut \\
\cline{2-5}

& \parbox{3cm} {$SmallGroup(96, 187)=(C_{2}\times S_{4}):C_{2}$} & $A_{4}\cong Group([ f4, f5, f6 ])$ & 8 & 3\tstrut\bstrut \\
\cline{2-5}

& \parbox{3cm} {$SmallGroup(96, 194)=C_{2}\times (A_{4}:C_{4})$} & $A_{4}\cong Group([ f4, f5, f6 ])$ & 8 & 3 \tstrut\bstrut \\
\cline{2-5}

& \parbox{3.5cm}{$SmallGroup(96, 195)=(C_{2}\times C_{2} \times A_{4}):C_{2}$} & $A_{4}\cong Group([ f4, f5, f6 ])$ & 8 & 3\tstrut\bstrut \\
\cline{2-5}

& \parbox{3cm} {$SmallGroup(96, 226)=(C_{2}\times C_{2} \times S_{4})$} & $A_{4}\cong Group([ f4, f5, f6 ])$ & 8 & 3 \tstrut\bstrut \\
\cline{2-5}

& \parbox{3cm} {$SmallGroup(272, 50)=C_{17}:C_{16}$} & $D_{34}\cong Group([ f4, f5 ])$ & 17 & 16 \tstrut\bstrut \\
\cline{2-5}

& \multirow{3}{*}{\parbox{3.5cm}{$SmallGroup(144, 120)=((C_{3}\times C_{3}):C_{4}):C_{4}$}} & \parbox{3.5cm}{$(C_{3}\times C_{3}):C_{2}\cong Group([ f3*f4, f5, f6 ])$} & 9 & 8 \tstrut \\
&  & \parbox{3.5cm}{$(C_{3}\times C_{3}):C_{2}\cong Group([ f3, f5, f6 ])$} & 9 & 8\bstrut \\
\cline{2-5}

& \parbox{3.5cm}{$SmallGroup(144, 182)=((C_{3}\times C_{3}):C_{8}):C_{2}$} & \parbox{3.5cm}{$(C_{3}\times C_{3}):C_{2}\cong Group([ f4, f5, f6 ])$} & 9 & 8 \tstrut\bstrut \\
\cline{2-5}

& \multirow{2}{*}{\parbox{3.5cm}{$SmallGroup(144, 185)=C_{2}\times ((C_{3}\times C_{3}):C_{8})$}} & \parbox{3.5cm}{$(C_{3}\times C_{3}):C_{2}\cong Group([ f2*f4, f5, f6 ])$} & 9 & 8\tstrut\\
&  & \parbox{3.5cm}{$(C_{3}\times C_{3}):C_{2}\cong Group([ f4, f5, f6 ])$} & 9 & 8\bstrut \\
\cline{2-5}

& \multirow{2}{*}{\parbox{3.5cm}{$SmallGroup(144, 187)=C_{2}\times ((C_{3}\times C_{3}):Q_{8})$}} & \parbox{3.5cm}{$(C_{3}\times C_{3}):C_{2}\cong Group([ f3*f4, f5, f6 ])$} & 9 & 8 \tstrut \\
&  & \parbox{3.5cm}{$(C_{3}\times C_{3}):C_{2}\cong Group([ f4, f5, f6 ])$} & 9 & 8\bstrut \\
\cline{2-5}

& \multirow{2}{*}{\parbox{3cm} {$SmallGroup(80, 28)=(C_{5}:C_{8}):C_{2}$}} & $D_{10}\cong Group([ f2*f3*f4, f5 ])$ & 5 & 4 \tstrut\\
&  & $D_{10}\cong Group([ f2*f3, f5 ])$ & 5 & 4\bstrut \\
\cline{2-5}

& \multirow{2}{*}{\parbox{3cm} {$SmallGroup(80, 30)=C_{4}\times (C_{5}:C_{4})$}} & $D_{10}\cong Group([ f3*f4, f5 ])$ & 5 & 4 \tstrut\\
&  & $D_{10}\cong Group([ f3, f5 ])$ & 5 & 4\bstrut \\
\cline{2-5}

& \multirow{2}{*}{\parbox{3cm} {$SmallGroup(80, 31)=C_{20}:C_{4}$}} & $D_{10}\cong Group([ f3*f4, f5 ])$ & 5 & 4\tstrut\\
&  & $D_{10}\cong Group([ f3, f5 ])$ & 5 & 4\bstrut \\
\cline{2-5}

& \multirow{2}{*}{\parbox{3.5cm}{$SmallGroup(80, 34)=(C_{2}\times (C_{5}:C_{4})):C_{2}$}} & $D_{10}\cong Group([ f3*f4, f5 ])$ & 5 & 4\tstrut\\
&  & $D_{10}\cong Group([ f3, f5 ])$ & 5 & 4\bstrut \\
\cline{2-5}

& \multirow{4}{*}{\parbox{3.5cm}{$SmallGroup(80, 50)=C_{2}\times C_{2}\times (C_{5}:C_{4})$}} & $D_{10}\cong Group([ f2*f3*f4, f5 ])$ & 5 & 4\tstrut\\
&  & $D_{10}\cong Group([ f2*f4, f5 ])$ & 5 & 4\\
&  & $D_{10}\cong Group([ f3*f4, f5 ])$ & 5 & 4\\
&  & $D_{10}\cong Group([ f4, f5 ])$ & 5 & 4\bstrut \\
\cline{2-5}

& \multirow{2}{*}{\parbox{3cm} {$SmallGroup(48, 4)=C_{8}\times S_{3}$}} & $S_{3}\cong Group([ f1*f4, f5 ])$ & 3 & 2\tstrut\\
&  & $S_{3}\cong Group([ f1, f5 ])$ & 3 & 2\bstrut \\
\cline{2-5}

& \multirow{4}{*}{\parbox{3cm} {$SmallGroup(48, 35)=C_{2}\times C_{4} \times S_{3}$}} & $S_{3}\cong Group([ f1*f2*f4, f5 ])$ & 3 & 2 \tstrut\\
&  & $S_{3}\cong Group([ f1*f2, f5 ])$ & 3 & 2\\
&  & $S_{3}\cong Group([f1*f4, f5 ]$ & 3 & 2\\
&  & $S_{3}\cong Group([ f1, f5 ])$ & 3 & 2\bstrut \\
\cline{2-5}

& \multirow{2}{*}{\parbox{3cm} {$SmallGroup(48, 38)=D_{8}\times S_{3}$}} & $S_{3}\cong Group([ f1*f4, f5 ])$ & 3 & 2 \tstrut \\
&  & $S_{3}\cong Group([ f1, f5 ])$ & 3 & 2\bstrut \\
\cline{2-5}

& \multirow{2}{*}{\parbox{3cm} {$SmallGroup(48, 40)=Q_{8}\times S_{3}$}} & $S_{3}\cong Group([ f1*f4, f5 ])$ & 3 & 2 \tstrut\\
&  & $S_{3}\cong Group([ f1, f5 ])$ & 3 & 2\bstrut \\
\cline{2-5}

& \multirow{8}{*}{\parbox{3cm} {$SmallGroup(48, 51)=C_{2}\times C_{2}\times C_{2}\times S_{3}$}} & $S_{3}\cong Group([ f1*f2*f3*f4, f5 ])$ & 3 & 2\tstrut\\
&  & $S_{3}\cong Group([ f1*f2*f3, f5 ])$ & 3 & 2\\
&  & $S_{3}\cong Group([ f1*f2*f4, f5 ])$ & 3 & 2\\
&  & $S_{3}\cong Group([ f1*f2, f5 ])$ & 3 & 2\\
&  & $S_{3}\cong Group([f1*f3*f4, f5 ]$ & 3 & 2\\
&  & $S_{3}\cong Group([ f1*f3, f5 ])$ & 3 & 2\\
&  & $S_{3}\cong Group([f1*f4, f5 ]$ & 3 & 2\\
&  & $S_{3}\cong Group([ f1, f5 ])$ & 3 & 2\bstrut \\
\hline
\textbf{9} & \parbox{3cm} {$SmallGroup(342, 7)=(C_{19}:C_{9}):C_{2}$} & $D_{38}\cong Group([ f1, f4 ])$ & 19 & 18 \tstrut\bstrut \\
\cline{2-5}
& \parbox{3cm} {$SmallGroup(126, 1)=(C_{7}:C_{9}):C_{2}$} & $D_{14}\cong Group([ f1, f4 ])$ & 7 & 6 \tstrut\bstrut \\
\cline{2-5}
& \parbox{3.5cm}{$SmallGroup(126, 7)=C_{3} \times ((C_{7}:C_{3}):C_{2})$} & $D_{14}\cong Group([ f1, f4 ])$ & 7 & 6 \tstrut\bstrut \\
\cline{2-5}
& \parbox{3cm} {$SmallGroup(54, 4)=C_{9} \times S_{3}$} & $S_{3}\cong Group([ f1, f4 ])$ & 3 & 2 \tstrut\bstrut \\
\cline{2-5}
& \parbox{3cm} {$SmallGroup(54, 12)=C_{3} \times C_{3} \times S_{3}$} & $S_{3}\cong Group([ f1, f4 ])$ & 3 & 2 \tstrut\bstrut \\
\hline
\textbf{10} & \parbox{3cm} {$SmallGroup(120, 37)=C_{5}\times S_{4}$} & $A_{4}\cong Group([ f3, f4, f5 ])$ & 8 & 3 \tstrut\bstrut \\
\cline{2-5}
& \parbox{3cm} {$SmallGroup(120, 38)=(C_{5}\times A_{4}):C_{2}$} & $A_{4}\cong Group([ f2, f4, f5 ])$ & 8 & 3 \tstrut\bstrut \\
\cline{2-5}
& \multirow{2}{*}{\parbox{3.5cm}{$SmallGroup(220, 7)=C_{2}\times ((C_{11}:C_{5}):C_{2})$}} & $D_{22}\cong Group([ f1*f2, f4 ])$ & 11 & 10 \tstrut\\
&  & $D_{22}\cong Group([ f1, f4 ])$ & 3 & 2\bstrut \\
\cline{2-5}
& \parbox{3cm} {$SmallGroup(100, 9)=C_{5} \times (C_{5}:C_{4})$} & $D_{10}\cong Group([ f3,f4 ])$ & 5 & 4 \tstrut\bstrut \\
\cline{2-5}
& \parbox{3cm} {$SmallGroup(100, 10)=(C_{5} \times C_{5}):C_{4}$} & $D_{10}\cong Group([ f2,f4 ])$ & 5 & 4\tstrut\bstrut \\
\cline{2-5}
& \parbox{3cm} {$SmallGroup(60, 8)=S_{3}\times D_{10}$} & $S_{3}\cong Group([ f2,f3 ])$ & 3 & 2 \tstrut\bstrut \\
\cline{2-5}
& \multirow{2}{*}{\parbox{3cm} {$SmallGroup(60, 11)=C_{10}\times S_{3}$}} & $S_{3}\cong Group([ f1*f2, f4 ])$ & 3 & 2 \tstrut\\
&  & $S_{3}\cong  Group([ f1, f4 ])$ & 3 & 2\bstrut \\
\hline
\textbf{11}& \parbox{3cm} {$SmallGroup(506, 1)=(C_{23}:C_{11}):C_{2}$} & $D_{46}\cong Group( [ f1, f3 ])$ & 23 & 22 \tstrut\bstrut \\
\cline{2-5}
& \parbox{3cm} {$SmallGroup(66, 1)=C_{11}\times S_{3}$} & $S_{3}\cong Group([ f1, f3 ])$ & 3 & 2 \tstrut\bstrut \\
\hline
\textbf{13}& \parbox{3.8cm} {$SmallGroup(702, 47)=((C_{3}\times C_{3}\times C_{3}):C_{13}):C_{2}$} & \parbox{3.5cm}{$(C_{3}\times C_{3}\times C_{3}):C_{2}\cong Group([ f1, f3, f4, f5 ])$} & 27 & 26\tstrut\bstrut \\
\cline{2-5}
& \parbox{3cm} {$SmallGroup(78, 3)=C_{13}\times S_{3}$} & $S_{3}\cong Group([ f1, f3 ])$ & 3 & 2\tstrut\bstrut \\
\hline
\textbf{17}& \parbox{3cm} {$SmallGroup(102, 1)=C_{17}\times S_{3}$} & $S_{3}\cong Group([ f1, f3 ])$ & 3 & 2 \tstrut\bstrut \\
\hline
\end{longtable}
\end{center}
\end{footnotesize}

\FloatBarrier

{\bf Remark 2.} For every index, the symmetric group $S_{3}$ always appears as a normal subgroup. This happens because we can always construct groups in the following way: $G=N\times A$ where $N\unlhd G$ and $A$ is abelian. Consequently, the $G$-classes of $N$ are exactly the classes of $N$. Thus, by taking $N \cong S_{3}$, the $G$-classes of $N$ are exactly the classes of $S_{3}$, which sizes are $\lbrace 1,2,3\rbrace$. In fact, for $n=7$ and $n=17$ this is the unique possibility. \\

\section{Landau's extension for prime-power order elements}

Let $G$ be a finite group and let $N \unlhd G$. We denote by $k(G)$ the number of conjugacy classes of $G$ and by $k_{G}(N)$ the number of $G$-conjugacy classes of elements of $N$. Analogously, we denote by $kpp(G)$ the number of conjugacy classes of prime-power order elements of $G$ and by $kpp_{G}(N)$ the number of $G$-classes of prime-power order elements of $N$.\\

The aim of this section is to extend Landau's result for elements of prime-power order lying in a normal subgroup of a finite group. However, if we restrict our attention to just non-central $G$-classes  of prime-power order elements contained in a normal subgroup $N$ of $G$, one can easily see that $|N|$ and $|G|$ cannot be upper bounded in terms of the number of such classes  although the index $|G:N|$ is fixed. For example, suppose that $N\unlhd G$ with $|G:N|=n$ and such that $N$ has just one non-central $G$-class (necessarily of prime-power order elements), as in any example of Section 3.1. We know then that $N$ is a $p$-group for some prime $p$ and we can take an arbitrary abelian finite $p'$-group $H$ and construct $N_0=N\times H$ and $G=G\times H$. It clearly follows that $N_0\unlhd G_0$, with index $n$ too,  and $N_0$ contains exactly one non-central $G_0$-class of prime-power order elements. However, $|N_0|$ and $|G_0|$ may be freely large and consequently, we can claim that there are infinitely many finite groups $G$ having a normal subgroup of fixed index $n$ and containing a  fixed number of non-central $G$-conjugacy classes of prime-power order elements.\\

Therefore, in this section, we will consider all $G$-classes, central and non-central, contained in a normal subgroup, and will prove that there are finitely many groups $G$ having a normal subgroup $N$ of index $n$ and  with $kpp_G(N)=k$. We also provide an upper bound for $|G|$ and $|N|$ by means of $k$ in the case in which $G$ is solvable.\\

As we have pointed in the Introduction, H\'ethelyi and K\"ulshammer prove in \cite{Kulshammer} the following extension of Landau's theorem for prime-power order elements. Nevertheless, they did not give any numerical expression for the upper bound. \\

{\bf Theorem 4.1}. {\it For a positive integer $k$, there are only finitely many finite groups, up to isomorphism, satisfying that $k=kpp(G)$.}\\

{\it Proof}. This is Theorem 1.1 of \cite{Kulshammer}. $\Box$ \\

The above theorem is proved as a direct consequence of the following, which uses the classification of the Finite Simple Groups. \\

{\bf Lemma 4.2}. {\it There exists a function $\alpha:\mathbb{N}\longrightarrow \mathbb{N}$ with the following property: whenever $k$ is a positive integer and $G$ is a finite group with $kpp(G)=k$ then $|G|\leqslant \alpha(k)$}.\\

{\it Proof}. This is Lemma 1.5 of \cite{Kulshammer}. $\Box$ \\

Analogous results can be obtained when we work with conjugacy classes lying in a normal subgroup.
For this purpose, we observe the following known elementary facts. \\

{\bf Lemma 4.3}. {\it Let $G$ a finite group and let $N\unlhd G$. For every $x \in N$, the class $x^{G}$ decomposes into the union of exactly $|G:N{\rm \textbf{C}}_{G}(x)|$ conjugacy classes of $N$}.\\

{\it Proof}. By taking a set of representatives of the right coclasses of $N$ in $G$, we can write $x^{G}=\cup_{i=1}^{k}x_{i}^{N}$ for certain $x_i\in N$, which are $G$-conjugate. Furthermore, the fact that the $x_i$ are $G$-conjugate implies that the classes $x_i^N$ all have the same cardinality, and hence $|x^{G}|=k|x^{N}|$. On the other hand,
\begin{equation*}
|x^{G}|=|G:{\rm \textbf{C}}_{G}(x)N||{\textbf{C}}_{G}(x)N:{\rm \textbf{C}}_{G}(x)|=|G:{\rm \textbf{C}}_{G}(x)N||x^{N}|,
\end{equation*}
so we get $|G:{\rm \textbf{C}}_{G}(x)N|=k$, as wanted. $\Box$\\

{\bf Theorem 4.4}.  {\it Let $G$ a finite group and $N\unlhd G$. Then}
\begin{enumerate}
\item[(i)]  {\it$ k_{G}(N)\geq k(N)/|G:N|$,}
\item[(ii)]  {\it $kpp_{G}(N)\geq kpp(N)/|G:N|$.}
\end{enumerate}

{\it Proof}. (i) By Lemma 4.3, each $G$-conjugacy class of an element $x\in N$ decomposes in exactly $|G:N{\rm \textbf{C}}_{G}(x)|$ classes of $N$.  Since this number divides $|G:N|$,  we deduce that every $G$-class of every element in $N$ decomposes in at most in $|G:N|$ classes of $N$. Moreover, the $N$-classes appearing in such a decomposition are trivially distinct, and hence, the inequality follows.

(ii) The same argument of i) works, just by restricting to prime-power order elements of $N$. $\Box$\\

Now, we are ready to prove the first part of Theorem B of the Introduction, which is Theorem 4.6.\\

{\bf Lemma 4.5}. {\it There exists a function $\gamma : \mathbb{N}\times \mathbb{N} \longrightarrow \mathbb{N}$ with the following property. If $n, l \in \mathbb{N}^{+}$ and $G$ is a finite group with $N\unlhd G$ such that $|G:N|=n$ and $kpp_{G}(N)=l$, then $|G|\leqslant \gamma(n,l)$}.\\

{\it Proof}. By Theorem 4.4  we have $nl=kpp_{G}(N)|G:N|\geq kpp(N)$. Now,  for every $i$, $1\leqslant i \leqslant nl$, if we assume that $kpp(N)=i$ , by Lemma $4.2$, we know that there exists $\alpha: \mathbb{N}\longrightarrow \mathbb{N} $ such that $|N|\leqslant \alpha(i)$. Therefore,  if we define  $\gamma(n,l)=|G:N| {\rm max}\{\alpha(1), \cdots, \alpha(nl)\}$, we conclude that $|G|\leqslant \gamma(n,l)$. $\Box$\\

{\bf Theorem 4.6}. {\it For any two positive integers $k, n$, there are only finitely many finite groups, up to isomorphism, with a normal subgroup of index $n$ satisfying that $k=kpp_{G}(N)$.}\\

{\it Proof}. This is a direct consequence of Lemma 4.5. $\Box$\\

 From now on we will assume that $G$ is a solvable group and we will get first an specific upper bound for $|G|$ when $kpp(G)=k$. The idea of our proof is inspired by Lemma 1.5 of \cite{Kulshammer}. Afterwards, if we take a normal subgroup $N$ of $G$ with a fixed index and with $kpp_G(N)=k$, we will obtain in Corollary 4.9 a formula for bounding $|N|$ and $|G|$, so Theorem B will be proved.\\

 In order to get the upper bound for $|G|$, we need one more result of \cite{Kulshammer} concerning to $G$-classes and factor groups. This holds for any finite group, not only for solvable groups.\\

{\bf Lemma 4.7}. {Let $N$ be a normal subgroup of a finite group $G$. Then $kpp(G/N)< kpp(G)$ unless $N=1.$}\\

{\it Proof}. This is Lemma 1.2(ii) of \cite{Kulshammer}. $\Box$\\

{\bf Theorem 4.8}. {\it Let G be a finite solvable group such that $kpp(G)=k$. Then $$|G|\leqslant \gamma(k)$$ where $\gamma$ is defined as follows: $\gamma(1)=1$ and $\gamma(k)=k\gamma(k-1)^{2}$  for every $k\geqslant2$. Consequently, $$|G|\leqslant \prod_{i=0}^{k-1}(k-i)^{2^{i}}.$$ }

{\it Proof}. We argue by induction on $k$ to prove both that $|G|\leqslant \gamma(k)$ and that $\gamma$ is an increasing function in $[1, k]$. When $k=1$ it is trivial that $|G|=1$, so we define $\gamma(1)=1$. If $k=2$, we take $N$ a minimal normal subgroup $G$, and by applying Lemma 4.7, since $N\neq1$, we have $kpp(G/N)< kpp(G)=2$. Then,  $kpp(G/N)=1$, that is, $N=G$. On the other hand, since $N$ is abelian, we certainly have $|N|=2$. Therefore, $|G|\leqslant 2= 2\gamma(1)^2=:\gamma(2)$. Also, in particular we have $\gamma(1)<\gamma(2)$. Now, assume by induction that when $kpp(G)=k$ then $|G|\leqslant k\gamma(k-1)^2$ and that $\gamma$ is an increasing function in $[1, k]$. Suppose that $kpp(G)=k+1$ and let $N$ be a minimal normal subgroup of $G$.  Again by Lemma 4.7, we have $kpp(G/N)< kpp(G)=k+1$. Then, $|G/N|\leqslant {\rm max}\{\gamma(1),\ldots , \gamma(k)\}=\gamma(k)$ by the inductive hypothesis. Furthermore, $N$ trivially contains at most $k+1$ $G$-conjugacy classes of elements of prime-power order. Each of these splits into at most $|G/N|$ $N$-classes of elements of prime-power order. Thus, $N$ has at most $(k+1)\gamma(k)$ conjugacy classes of elements of primer power order by Theorem 4.4(ii). As $N$ is an abelian prime-power order group, we conclude that $|N|\leqslant (k+1)\gamma(k)$. Therefore, $|G|=|N||G/N| \leqslant (k+1)\gamma(k)^2=:\gamma(k+1)$ and also, in particular, $\gamma(k)<\gamma(k+1)$. The latter inequality of the statement follows easily. $\Box$ \\

{\bf Corollary 4.9}. {\it Let G be a finite solvable group and let $N\trianglelefteq G$ such that $|G:N|=n$ and $kpp_{G}(N)=k$.  Then $$|N|\leqslant \prod_{i=0}^{nk-1}(nk-i)^{2^{i}}$$ and $$|G|\leqslant n\prod_{i=0}^{nk-1}(nk-i)^{2^{i}}.$$}

{\it Proof}. By Theorem 4.4, we know that $kpp_{G}(N)\geq kpp(N)/|G:N|$, so $kpp(N)\leqslant |G:N|kpp_{G}(N)=nk$. Thus, by applying Theorem 4.8 to $N$ we get $|N|\leqslant \gamma(nk)=\prod_{i=0}^{nk-1}(nk-i)^{2^{i}}$. This implies that $|G|\leqslant n|N|=n\gamma(nk)\leqslant n\prod_{i=0}^{nk-1}(nk-i)^{2^{i}}.$ $\Box$\\

We would like to remark that when $G$ is solvable, M. Lewis proves in \cite{Lewis}, as a consequence of Theorem 4.1, that the number of conjugacy classes of elements of prime-power order of $G$ is less than or equal to the number of irreducible characters of $G$  with values in fields where $\mathbb{Q}$ is extended by prime-power roots of the unity.\\

We finish this section by showing in Table 5 a classification, obtained with GAP, of those solvable groups $G$  with $kpp(G)=k$ for small values of $k$. We note that for $k=5$ the bound of Corollary 4.9 is too large to be used with GAP, but we may improve it by using the same argument as in Theorem 4.8. Let $N$ be a minimal normal subgroup of $G$. By Lemma 4.7 we know that $kpp(G/N)<kpp(G)=5$, so $kpp(G/N)\leq 4$. We observe in Table 5 that the largest size of $|G|$ for $k=4$ is $12$. Consequently, $|G/N|\leq 12$. Furthermore, $N$ trivially contains at most $5$ $G$-classes of elements of primer-power order and each of these splits into at most $|G/N|$ $N$-classes. Thus, $N$ has at most $5\times 12=60$ conjugacy classes of elements of prime power order by Theorem 4.4(ii). As a consequence, $|G|=|N||G/N|\leq60\times 12=720$. This bound is computationally appropriated.

\FloatBarrier
\begin{center}
\begin{longtable}[h]{l l }
\caption{Classification of solvable groups with small $kpp(G)$.}\\
\hline
 $k$ & $G$ \\
\hline \hline
\textbf{2} & $SmallGroup(2,1)=C_{2}$\tstrut\bstrut\\
\hline
\textbf{3} & $SmallGroup(3,1)=C_{3}$, $SmallGroup(6,1)=S_{3}$  \tstrut\bstrut\\
\hline
\textbf{4}&  $SmallGroup(4,1)=C_{4}$, $SmallGroup(4,2)=C_{2}\times C_{2}$,\tstrut\\ &$SmallGroup(6,2)=C_{6}$, $SmallGroup(10,1)=D_{10}$,\tstrut\\ &$SmallGroup(12,3)=A_{4}$\tstrut\bstrut\\
\hline
\textbf{5}&  $SmallGroup(5,1)=C_{5}$, $SmallGroup(8,3)=Q_{8}$, \tstrut\\ &$SmallGroup(8,4)=D_{8}$, $SmallGroup(12,1)=C_{3}:C_{4}$,\tstrut\\ &$SmallGroup(12,4)=D_{12}$, $SmallGroup(14,1)=D_{14}$, \tstrut\\ &$SmallGroup(20,3)=C_{5}:C_{4}$, $SmallGroup(21,1)=C_{7}:C_{3}$,\tstrut\\
&$SmallGroup(24,3)=SL(2,3)$, $SmallGroup(24,12)=S_{4}$, \tstrut\\ &$SmallGroup(30,3)=D_{30}$, $SmallGroup(42,1)=(C_7 : C_3) : C_2$\tstrut\bstrut\\
\hline
\end{longtable}
\end{center}

\FloatBarrier

\end{document}